\documentclass{article}

\usepackage{amsthm}
\usepackage{amsmath}
\usepackage{amsfonts}
\usepackage{amssymb}
\usepackage{calc}

\usepackage[all,ps,dvips,arc,knot,frame]{xy}     

\xyresetcatcodes

\def \R{{\mathbb R}}

\def \Z{{\mathbb Z}}
\def \A{{\mathcal A}}

\def \B{{\mathcal B}}
\def \V{{\mathcal V}}

\def \bracket#1{\bigl< #1 \bigr>}

\def \bf{\bfseries}

\def \ur#1{\,\overline{\,#1\vphantom{ly}\,}\raisebox{-.81ex}{\rule{.4pt}{2.8ex}}\,}
\def \lr#1{\,\underline{\,#1\vphantom{ly}\,}\raisebox{-.81ex}{\rule{.4pt}{2.8ex}}\,}
\def \ul#1{\,\raisebox{-.81ex}{\rule{.4pt}{2.8ex}}\overline{\ \vphantom{ly}#1\,}\,}
\def \ll#1{\,\raisebox{-.81ex}{\rule{.4pt}{2.8ex}}\underline{\ \vphantom{ly}#1\,}\,}
\def \Bur#1{\,\overline{\,#1\vphantom{\raisebox{-.90ex}{\rule{.4pt}{2.94ex}}}\,}\raisebox{-1.28ex}{\rule{.4pt}{3.72ex}}\,}
\def \Blr#1{\,\underline{\,#1\vphantom{\raisebox{-.90ex}{\rule{.4pt}{2.94ex}}}\,}\raisebox{-1.28ex}{\rule{.4pt}{3.72ex}}\,}
\def \Bul#1{\,\raisebox{-1.28ex}{\rule{.4pt}{3.72ex}}\overline{\ \vphantom{\raisebox{-.90ex}{\rule{.4pt}{2.94ex}}}#1\,}\,}
\def \Bll#1{\,\raisebox{-1.28ex}{\rule{.4pt}{3.72ex}}\underline{\ \vphantom{\raisebox{-.90ex}{\rule{.4pt}{2.94ex}}}#1\,}\,}
\def \Bbur#1{\,\overline{\,#1\vphantom{\raisebox{-1.37ex}{\rule{.4pt}{3.85ex}}}\,}\raisebox{-1.75ex}{\rule{.4pt}{4.6ex}}\,}
\def \Bblr#1{\,\underline{\,#1\vphantom{\raisebox{-1.37ex}{\rule{.4pt}{3.85ex}}}\,}\raisebox{-1.75ex}{\rule{.4pt}{4.6ex}}\,
}

\theoremstyle{plain}
\newtheorem{thm}{Theorem}[section]
\newtheorem{cor}[thm]{Corollary}
\newtheorem{lemma}[thm]{Lemma}

\theoremstyle{definition}
\newtheorem{defn}[thm]{Definition}
\theoremstyle{remark}

\begin{document}

\date{}

\title{\Large\bf Biquandles for Virtual Knots }
\author{
David Hrencecin \\
Department of Mathematics, Statistics \\
and Computer Science (m/c 249)    \\
851 South Morgan Street   \\
University of Illinois at Chicago\\
Chicago, Illinois 60607-7045\\
$<$djhren at gmail dot com$>$\\
\\
and \\
\\
Louis H. Kauffman\\
Department of Mathematics, Statistics \\
and Computer Science (m/c 249)    \\
851 South Morgan Street   \\
University of Illinois at Chicago\\
Chicago, Illinois 60607-7045\\
$<$kauffman at uic dot edu$>$
}

\maketitle

\label{C:Summary}

Virtual knot theory can be viewed as the theory of abstract Gauss
codes. A Gauss code for a classical knot is obtained by walking
along the knot diagram and recording the names and ``states" of
the crossings in the order that the knot is traversed. This
abstract information can be used as a substitute for the knot.
Arbitrary Gauss codes are not necessarily planar, but can be used
to generalize knots.  That generalization is virtual knot theory.
We recommend~\cite{LouVirt} as an introduction for the reader who
is not already familiar with virtual knot theory.

This paper studies an algebraic invariant of virtual knots called
the biquandle.  The biquandle generalizes the fundamental group and
the quandle of virtual knots. The approach taken in this paper to
the biquandle emphasizes understanding its structure in terms of
compositions of morphisms, where elementary morphisms are associated
to oriented classical and virtual crossings in the diagram.

Section~\ref{C:Definitions} gives definitions of virtual knots and
reminds the reader of the basic result that classical knot theory
embeds in virtual knot theory. Section~\ref{C:Biquandle} defines the
biquandle. Section~\ref{S:SimplerBQ} shows how to associate
morphisms to the crossings. Section~\ref{S:ReversingBQ} gives our
basic formulae for inverting the morphism corresponding to a
crossing. We use this method to show that the biquandle descriptions
for the virtual Hopf link and its mirror image are identical.  Of
course, this shows that biquandles do not classify virtual links. We
also use this method to show that the biquandle is invariant under
AD inversion of the diagram (see Definition~\ref{D:ADinv}). This
method of inverting a diagram is defined in the paper.
Section~\ref{S:BQsViaBraids} explains how to use braids to compute
the biquandle. We also prove that $BQ(K)\cong BQ(K^\uparrow)$. The
vertical mirror image $K^\uparrow$ is one of various kinds of mirror
images that can be defined for virtual knots.  The section ends with
an example (see Figure~\ref{Fi:Mirror}) showing the difference
between two such types of mirror images. Section~\ref{S:ABQ} applies
our methods to calculations of the Alexander biquandle.
Section~\ref{C:BQKishino} discusses the structure of the biquandle
description of the Kishino knot, and uses our framework to give a
proof of its non-triviality.

\section{Definitions} \label{C:Definitions}

We define a knot $K$ in the combinatorial sense, as a class of
diagrams which represent a generic projection of an embedding $S^1\to S^3$
or $S^1\sqcup\dots\sqcup S^1\to S^3$.  Each circle represented in such a diagram
is called a component, and if there is more than one component, the diagram (or related class)
is sometimes referred to as a link.
A strand in a knot diagram is the projection of a connected interval in the
embedded curve.
In a knot diagram,
Certain arcs of the projection of the curve embedded
in space are eliminated to form a knot diagram, creating broken strands
that indicate the over and under crossings.
At each crossing in the diagram there is an indicated over strand and a broken
under strand.
In this sense, there are two local strands at any crossing
in a diagram.  Whenever we refer to a strand of a crossing, we mean one
of these two local strands.

\begin{defn}
Two knot diagrams $K$ and $K'$ are said to be \emph{equivalent} when there is a
finite sequence of the following (Reidemeister) moves which transform $K$ into $K'$:
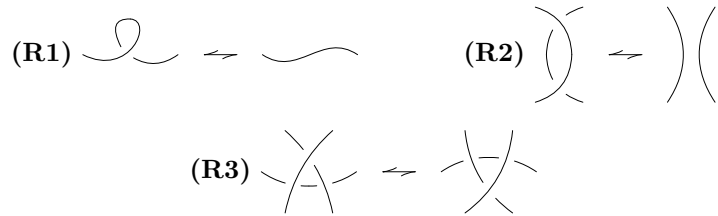
$$
\hbox{\bf(R1)\ }
\xy 0;<1pc,0pc>:
0;(1.5,1.125)**\crv{(1.5,-0.75)&(2.25,1.075)}
,(1.5,1.125);(1.2,0.25)**\crv{(0.750,1.15)}
,(1.6,-0.1);(3,0)**\crv{(2.3,-0.50)}
\endxy \quad \xy 0;<1pc,0pc>:
\ar @{^{<}-_{>}} ,(0,0.25);(1,0.25)
\endxy \quad
\xy 0;<1pc,0pc>:*=dir{}
,0;(3,0)**\crv{(1.125,-0.75)&(1.875,0.75)}
\endxy \qquad\qquad
\hbox{\bf(R2)\ }
\xy 0;<1pc,0pc>:
,\ar @{-} @/_4pt/,(0,-1.5);(1.15,0)
,\ar @{-} @/_4pt/,(1.15,0);(0,1.5)
,\ar @{-} @/^4pt/|(0.45)\hole,(1.5,-1.5);(0.35,0)
,\ar @{-} @/^4pt/|(0.55)\hole,(0.35,0);(1.5,1.5)
\endxy \quad \xy 0;<1pc,0pc>:
\ar @{^{<}-_{>}} ,(0,0.25);(1,0.25)
\endxy \quad \xy 0;<1pc,0pc>:*=dir{}
\ar @{-} @/_/ ,(0,-1.5);(0,1.5)
\ar @{-} @/^/ ,(1.5,-1.5);(1.5,1.5)
\endxy
$$
$$
\hbox{\bf(R3)\ }
\xy 0;<1.5pc,0pc>:
a(0)="a1",a(60)="b1",a(120)="c1",a(180)="a2",a(240)="b2",a(300)="c2"
,\ar @{-}@/^1.5pt/|(0.69)\hole "a1";(0,-0.3)
,\ar @{-}@/^1.5pt/|(0.31)\hole (0,-0.3);"a2"
,\ar @{-}@/^3pt/ "b2";"b1"
,\ar @{-}@/^3pt/|(0.37)\hole "c1";"c2"
\endxy \quad \xy
0;<1pc,0pc>:\ar @{^{<}-_{>}} ,(0,0.25);(1,0.25)
\endxy \quad \xy 0;<1.5pc,0pc>:
a(0)="a1",a(60)="b1",a(120)="c1",a(180)="a2",a(240)="b2",a(300)="c2"
,\ar @{-}@/_1.5pt/|(0.65)\hole "a1";(0,0.3)
,\ar @{-}@/_1.5pt/|(0.35)\hole (0,0.3);"a2"
,\ar @{-}@/_4pt/ "b2";"b1"
,\ar @{-}@/_4pt/|(0.67)\hole "c1";"c2"
\endxy
$$
A \emph{knot} is an equivalence class of knot diagrams under the Reidemeister moves.
\end{defn}
A knot can have an orientation.  This is indicated on a knot diagram
by drawing an arrowhead on one or more strands in such a way that
each component has a consistent labeling. The Reidemeister moves for
an oriented knot are the same as for unoriented knots. We are free
to apply the moves without paying attention to the particular
orientations of the strands involved.

There is also the notion of more than one kind of Reidemeister
equivalence. \emph{Regular isotopy} is equivalence under only the R2
and R3 moves. \emph{Ambient isotopy} is equivalence under all three
moves. It is useful to distinguish between the two because there are
invariants only of regular isotopy.  For references on knot theory,
see ~\cite{LouOnKs},~\cite{Rolfson} and~\cite{Kawauchi}.

\begin{defn} \label{D:vk}
A \emph{virtual knot diagram} is a generic immersion $S^1\sqcup\dots\sqcup S^1\to \R^2$
such that each double point is labelled with either a real (over or under) or a virtual
(circled) crossing.  Thus a \emph{virtual knot} is a class of equivalent
virtual diagrams, where two virtual knot diagrams are said to be equivalent
when one can be transformed into the other by a finite sequence of real
Reidemeister moves (R1), (R2) and (R3), along with the following virtual moves:
$$
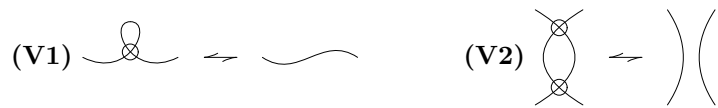
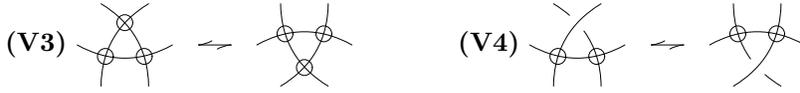
\hbox{\bf(V1)\ }
\xy 0;<1pc,0pc>:
0;(1.5,1.125)**\crv{(1.5,-0.75)&(2.25,1.125)} ?(0.45)*\cir<3pt>{}
,(1.5,1.125);(3,0)**\crv{(0.75,1.125)&(1.5,-0.75)}
\endxy \quad \xy 0;<1pc,0pc>:
\ar @{^{<}-_{>}} ,(0,0.25);(1,0.25)
\endxy \quad \xy
0;<1pc,0pc>:*=dir{}
,0;(3,0)**\crv{(1.125,-0.75)&(1.875,0.75)}
\endxy \qquad\qquad
\hbox{\bf(V2)\ }
\xy 0;<1pc,0pc>:
(0.75,0.94)*\cir<3pt>{}
,(0.75,-0.94)*\cir<3pt>{}
,\ar @{-} @/_15pt/,(0,-1.5);(0,1.5)
,\ar @{-} @/^15pt/,(1.5,-1.5);(1.5,1.5)
\endxy \quad \xy
0;<1pc,0pc>:
\ar @{^{<}-_{>}} ,(0,0.25);(1,0.25)
\endxy \quad \xy
0;<1pc,0pc>:*=dir{}
\ar @{-} @/_/ ,(0,-1.5);(0,1.5)
\ar @{-} @/^/ ,(1.5,-1.5);(1.5,1.5)
\endxy
$$
$$
\hbox{\bf(V3)\ }
\xy 0;<8.5pt,0pc>:
a(-30)*\cir<3pt>{},a(90)*\cir<3pt>{},a(210)*\cir<3pt>{}
,0;<1.5pc,0pc>:
a(0)="a1",a(60)="b1",a(120)="c1",a(180)="a2",a(240)="b2",a(300)="c2"
,\ar @{-}@/^5pt/,"a1";"a2"
,\ar @{-}@/^5pt/,"b2";"b1"
,\ar @{-}@/^5pt/,"c1";"c2"
\endxy \quad \xy
0;<1pc,0pc>:\ar @{^{<}-_{>}} ,(0,0.25);(1,0.25)
\endxy \quad \xy
0;<8.5pt,0pc>:
a(30)*\cir<3pt>{},a(150)*\cir<3pt>{},a(270)*\cir<3pt>{}
,0;<1.5pc,0pc>:
a(0)="a1",a(60)="b1",a(120)="c1",a(180)="a2",a(240)="b2",a(300)="c2"
,\ar @{-}@/_5pt/,"a1";"a2"
,\ar @{-}@/_5pt/,"b2";"b1"
,\ar @{-}@/_5pt/,"c1";"c2"
\endxy \qquad\qquad
\hbox{\bf(V4)\ }
\xy 0;<8.5pt,0pc>:
a(-30)*\cir<3pt>{},a(210)*\cir<3pt>{}
,0;<1.5pc,0pc>:
a(0)="a1",a(60)="b1",a(120)="c1",a(180)="a2",a(240)="b2",a(300)="c2"
,\ar @{-}@/^5pt/ "a1";"a2"
,\ar @{-}@/^5pt/ "b2";"b1"
,\ar @{-}@/^5pt/|(0.30)\hole "c1";"c2"
\endxy \quad \xy
0;<1pc,0pc>:\ar @{^{<}-_{>}} ,(0,0.25);(1,0.25)
\endxy \quad \xy
0;<8.5pt,0pc>:
a(30)*\cir<3pt>{},a(150)*\cir<3pt>{}
,0;<1.5pc,0pc>:
a(0)="a1",a(60)="b1",a(120)="c1",a(180)="a2",a(240)="b2",a(300)="c2"
,\ar @{-}@/_5pt/ "a1";"a2"
,\ar @{-}@/_5pt/ "b2";"b1"
,\ar @{-}@/_5pt/|(0.70)\hole "c1";"c2"
\endxy
$$
\end{defn}
As in the classical case, a virtual knot with multiple components is sometimes called
a virtual link.  For this paper, we will use the term `virtual knot' to refer to both
virtual knots and links.  When we do not wish to consider virtual crossings, we
will use the term `classical knot' to mean knots and links without virtual crossings.

An alternate definition to virtual Reidemeister equivalence is to allow the
classical Reidemeister moves with the addition of a more general ``detour move"
(See~\cite{LouVirt,LouVirt2}):
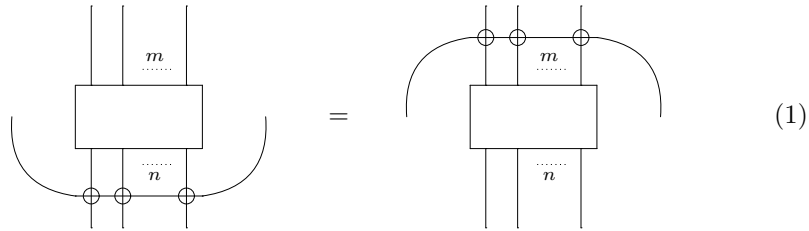
\begin{equation} \label{E:DetourMove}
\xy 0;/r1pc/:@={(0,1),(4,1),(4,-1),(0,-1)},s0="prev" @@{;"prev";**\dir{-}="prev"}
@i @={(.5,-2.5),(1.5,-2.5),(3.5,-2.5)} @@{*\cir<3pt>{}}
,\ar^{m}@{.} (2,1.5);(3,1.5),\ar_{n}@{.} (2,-1.5);(3,-1.5)
,\ar@{-} (.5,1);(.5,3.5),\ar@{-} (1.5,1);(1.5,3.5),\ar@{-} (3.5,1);(3.5,3.5)
,\ar@{-} (.5,-1);(.5,-3.5),\ar@{-} (1.5,-1);(1.5,-3.5),\ar@{-} (3.5,-1);(3.5,-3.5)
,\ar@{-} (0,-2.5);(4,-2.5),\ar@{-}@/_9pt/ (-2,0);(0,-2.5),\ar@{-}@/_9pt/ (4,-2.5);(6,0)
\endxy
\qquad = \qquad\quad
\xy 0;/r1pc/:@={(0,1),(4,1),(4,-1),(0,-1)},s0="prev" @@{;"prev";**\dir{-}="prev"}
@i @={(.5,2.5),(1.5,2.5),(3.5,2.5)} @@{*\cir<3pt>{}}
,\ar^{m}@{.} (2,1.5);(3,1.5),\ar_{n}@{.} (2,-1.5);(3,-1.5)
,\ar@{-} (.5,1);(.5,3.5),\ar@{-} (1.5,1);(1.5,3.5),\ar@{-} (3.5,1);(3.5,3.5)
,\ar@{-} (0,2.5);(4,2.5),\ar@{-}@/^9pt/ (-2,0);(0,2.5),\ar@{-}@/^9pt/ (4,2.5);(6,0)
,\ar@{-} (.5,-1);(.5,-3.5),\ar@{-} (1.5,-1);(1.5,-3.5),\ar@{-} (3.5,-1);(3.5,-3.5)
\endxy
\end{equation}
In the detour move, any number of strands may emanate from the top and bottom of the
tangle (represented by a box).  The idea is that in a virtual diagram, if we have an
arc with any number of consecutive virtual crossings, then we can cut that arc out
and replace it with another arc connecting the same points, provided that any crossings
on the new arc are also virtual.  It is easily seen that this yields the same equivalence.

Notice that there are some mixed moves which are not allowed.
Consider the following moves:
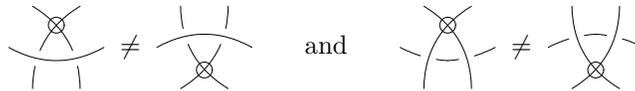
\[
\xy ,0;<1.5pc,0pc>:
<0pt,8.5pt>*\cir<3pt>{}
,a(0)="a1",a(60)="b1",a(120)="c1",a(180)="a2",a(240)="b2",a(300)="c2"
,\ar @{-}@/^5pt/"a1";"a2"
,\ar @{-}@/^5pt/|(0.30)\hole "b2";"b1"
,\ar @{-}@/^5pt/|(0.70)\hole "c1";"c2"
\endxy
\ \ne\
\xy ,0;<1.5pc,0pc>:
<0pt,-8.5pt>*\cir<3pt>{}
,a(0)="a1",a(60)="b1",a(120)="c1",a(180)="a2",a(240)="b2",a(300)="c2"
,\ar @{-}@/_5pt/"a1";"a2"
,\ar @{-}@/_5pt/|(0.70)\hole  "b2";"b1"
,\ar @{-}@/_5pt/|(0.30)\hole "c1";"c2"
\endxy
\qquad\hbox{and}\qquad
\xy ,0;<1.5pc,0pc>:
a(0)="a1",a(60)="b1",a(120)="c1",a(180)="a2",a(240)="b2",a(300)="c2"
,<0pt,8.5pt> *\cir<3pt>{},<0pt,-5pt>="a3"
,\ar @{-}@/^1pt/|(0.62)\hole "a1";"a3"
,\ar @{-}@/^1pt/|(0.39)\hole "a3";"a2"
,\ar @{-}@/^5pt/ "b2";"b1"
,\ar @{-}@/^5pt/ "c1";"c2"
\endxy
\ \ne\
\xy ,0;<1.5pc,0pc>:
a(0)="a1",a(60)="b1",a(120)="c1",a(180)="a2",a(240)="b2",a(300)="c2"
,<0pt,-8.5pt> *\cir<3pt>{},<0pt,5pt>="a3"
,\ar @{-}@/_1pt/|(0.62)\hole "a1";"a3"
,\ar @{-}@/_1pt/|(0.39)\hole "a3";"a2"
,\ar @{-}@/_5pt/ "b2";"b1"
,\ar @{-}@/_5pt/ "c1";"c2"
\endxy
\]
These are forbidden because they change the underlying Gauss code of a knot,
and the only changes we allow on the Gauss codes are the Reidemeister moves.

However, inclusion of the above left (overstrand) version
of this move has been studied in the form of welded braids~\cite{Fenn,FennEtal}
-- a generalization of braid theory preceding virtual knot theory.

There has been progress in applying well known knot invariants
to the virtual theory.  In~\cite{LouVirt}, Kauffman extended the
fundamental group, the Jones polynomial and classes of quantum link invariants
to virtual knots and gave examples of non-trivial virtual knots with trivial
Jones Polynomial and trivial\footnote{by trivial fundamental
group, we mean a group isomorphic to the integers} fundamental group.

One of the results we will assume is the following, proved by
Kauffman and by Goussarov, Polyak and Viro:
\begin{lemma}
If $K$ and $K'$ are classical knots which are equivalent under virtual Reidemeister
equivalence, then they are also equivalent under classical Reidemeister equivalence.
\end{lemma}

An immediate question that arises is how to determine when a virtual knot
is classical.  Many virtual knots are indeed non-classical, and in some cases
it not possible to determine this property using extensions of classical
invariants.
In this paper, we will discuss the biquandle, a recent invariant, which
turns out to be surprisingly useful in detecting non-classical knots.

\section{The Biquandle} \label{C:Biquandle}
The biquandle has a non-associative structure, and in order to
handle it, we will use a special operator notation system, whose
symbols are $\ur{\ }$, $\ul{\ }$, $\ll{\ }$ and $\lr{\ }.$  The
notation allows us to write expressions in a non-associative
setting without requiring lots of parentheses.  Suppose we have
an algebraic system with four binary operations,
$*$, $\bar{*}$, $\sharp$ and $\bar{\sharp}$.

If $X$ is any symbol-string representing an expression in the
algebra, and $Y$ is another such string, then the binary operations
on $X$ and $Y$ are expressed by
\begin{align*}
X\hbox{$*$}Y&=X\ur Y  &  X\sharp Y&=X\lr Y\\
X\bar{*}Y&=X\ul Y  &  X\bar{\sharp}Y&=X\ll Y \\
\end{align*}

It then follows that the associations between these symbols are expressed by:
\begin{align*}
(X*Y)*Z &= X\ur Y\ur Z \\
X*(Y*Z) &= X\Bur{Y\ur Z} \\
\intertext{ Use the same convention for all operators.  Here are
some mixed expressions: }
(X\,\bar{*}\,Y)*Z &= X\ul Y \ur Z \\
(X\sharp(Y\,\bar{*}\,(Z\bar{\sharp}W)))\bar{*}T &= X\Bblr{Y \Bul{Z\ll W}}\ul T \\
\end{align*}

In the operator formalism, the simplest expressions correspond to
the left-associated expressions in the binary algebra.  Note that when an
expression is associated differently, the over and under bars serve as
parentheses by extending over or under entire sub-expressions.

The biquandle~\cite{LouRad, Hren1, FennJK} is a generalization of
the quandle~\cite{FennRourke, Joyce}, which in turn is a generalization of the fundamental
group of a classical knot.

In the following, when we speak of an \emph{algebra}, we mean a set that is closed
under certain binary operations.  The rules specifying the properties of those
binary operations will be part of the definition of that algebra.

\begin{defn}[Biquandle Axioms]
A \emph{biquandle} $BQ$ is a non-as\-so\-ci\-a\-tive algebra with four binary operations
represented through the symbols
$\ur{\ }$, $\lr{\ }$, $\ul{\ }$ and $\ll{\ }$ which satisfy following axioms:
\begin{enumerate}
\item
For every element $a\in BQ$, there is an element $x\in BQ$ such that \ $a\ur x\lr a = a$. \\
The left/right variant of this statement must also be true, that is: \\
For every element $a\in BQ$, there is an element $x\in BQ$ such that \ $a\ul x\ll a = a$.
\item
For every element $a\in BQ$, there is an element $x\in BQ$ such that \ $x=a\ll x$ \ and \ $a=x\ul a$. \\
The left/right variant of this statement must also be true.
\item
For any $a,b\in BQ$, we have
\[
a\ =\ a\lr{b}\Bll{b\ur{a}}
\ =\  a\ur{b}\Bul{b\lr{a}}
\ =\  a\ul{b}\Bur{b\ll{a}}
\ =\  a\ll{b}\Blr{b\ul{a}}
\]
\item
For any $a,b\in BQ$, there is an element $x\in BQ$ such that
\[
x=a\Bur{b\ll x}, \qquad a=x\ul b \hbox{\qquad and \qquad} b=b\ll x\lr a
\]
The left/right variant of this statement must also be true.
\item
For any $a,b,c\in BQ$, the following equations hold, as well as the left/right variants:
\begin{align*}
a\ur b\ur c &= a\Bur{c\lr b}\Bur{b\ur c} \\
a\lr b\lr c &= a\Blr{c\ur b}\Blr{b\lr c} \\
a\lr b\Bbur{c\Blr{b\ur a}} &= a\ur c\Bblr{b\Bur{c\lr a}}
\end{align*}
\end{enumerate}
\end{defn}

We define a biquandle $BQ(K)$ associated with a knot $K$ as follows.
Start with an oriented knot diagram $K$.  Assign one generator to each
arc.  For each classical crossing, assign two relations according the
equations given in the diagrams below:
\begin{equation} \label{E:BQR}
\xy 0;<4pc,0pc>:
\ar ,(4,-0.5) *+!U{a\vphantom{l}};(5,0.5) *+!D{\quad d=a\lr{b}}
\ar|(0.5)\hole ,(5,-0.5) *+!U{b};(4,0.5) *+!D{c=b\ur{a}\quad}
\ar|(0.5)\hole ,(0,-0.5) *+!U{a\vphantom{l}};(1,0.5) *+!D{\quad d=a\ul{b}}
\ar ,(1,-0.5)*+!U{b} ;(0,0.5)*+!D{c=b\ll{a}\quad}
\endxy
\end{equation}
A diagram with $n$ real crossings gives rise to $2n$ generators and $2n$ relations.

An algebra $\A$ with four binary operations and $2n$ generators that satisfies the
biquandle axioms along with these $2n$ relations is said to be a biquandle for
the knot $K$.  In practice, it will often suffice to analyze the structure of any
hypothetical biquandle algebra that satisfies these generators and relations.

%

Our aim is to improve the description of the biquandle of a knot or link.
We will describe a simpler way to construct generators and relations.

\section{A Simplified Approach to Biquandle Computation} \label{S:SimplerBQ}
View each crossing (both real and virtual types) as a morphism on an ordered
pair of biquandle elements.  We will use the standard diagrammatic convention
that inputs are ordered in a counterclockwise direction, and outputs are
ordered clockwise.When we define a morphism related to a crossing, it will
act vertically up the page, mapping the inputs at the bottom of the crossing
to the outputs at the top.

Once we define the morphisms, we do not require that
the crossings be aligned upwards in the page.  Instead, we place
a mapping direction directly on the tangles within a diagram.  The direction
indicates input and output strands, as well as which morphisms to apply.
This will be made clear in the context of each diagram.

For example, let $\phi_u$ be the automorphism which describes how the crossing
$\raise1.3ex\xy 0;/r1pc/:\xoverv=<\endxy$ acts on $(a,b)$:
\begin{equation*}
\xy 0;/r2pc/:
    (-1,1)="a" *+!D{b\ur a} ,(1,1)="b" *+!D{a\lr b}
    ,(-1,-1)="c" *+!U{a} ,(1,-1)="d" *+!U{b}
    \xoverv~{"a"}{"b"}{"c"}{"d"}=<
\endxy
\text{\qquad defining $\phi_u$ via \qquad}
\begin{gathered}
(\phi_u)_1(a,b) = b\ur a \\
(\phi_u)_2(a,b) = a\lr b \\
\text{or}\\
\phi_u(a,b) = (b\ur a, a\lr b)
\end{gathered}
\end{equation*}
The subscript of $u$ indicates that the direction is upwards, or coinciding with
the direction of the strands in the crossing.

Since many knots contain multiple ``twists'' in succession, we view tangle
sums as acting on the two generators $(a,b)$ through compositions of the
elementary morphisms.  The following tangle becomes the
morphism $\phi_u^3$ when viewed as mapping from left to right:
\begin{equation} \label{E:BQcomp}
\xy 0;/r3pc/:0*+!R{a}
\htwist\htwist\htwist=>
,(0,-1) *+!R{b} ,(3,0) *+!L{(\phi_u^3)_1(a,b)}
,(3,-1) *+!L{(\phi_u^3)_2(a,b)}
\endxy
\end{equation}
This example is a special case of working with a braid.  Every virtual knot
is equivalent to the closure of some virtual braid, proved in
Theorem~\ref{L:VirtualAlexander}, a generalization of Alexander's theorem
for classical knots and links.

This point of view helps in many ways:
\begin{enumerate}
\item It decreases the number of generators and relations in a presentation
    by taking compositions of the morphisms related to a string of crossings.
\item It suggests that we might also consider the inverse automorphism
    associated with each crossing.
\item In examples such as the Alexander biquandle (see Section~\ref{S:ABQ}) where the crossings are linear maps, we
    may be able to define morphisms for every possible crossing rotation.
\end{enumerate}

In this paper,
a presentation of a biquandle in terms of generators and relations is
a description of some biquandle generated by these generators, which satisfies the additional
equations given by the crossings.  To discriminate such a structure from
from a universal algebraic presentation, we shall refer to it as a \emph{description} of
a biquandle.

In this paper, we have not discussed the universal algebra for biquandles.
In order to define a free biquandle, or a biquandle with generators and relations,
one needs to modify the axioms to accommodate the logic of this level of construction.
The details of this will be taken up elsewhere.

\section{Inverting the Biquandle Automorphisms} \label{S:ReversingBQ}
There are inverse equations for the biquandle.

\begin{thm} \label{T:BQ_inverse}
When assigning equations for the biquandle from any crossing in a virtual
diagram, we may instead use the following ``inverse equations'':
\begin{equation} \label{E:InvABQR}
\xy 0;<4pc,0pc>:
\ar ,(3,-0.5) *+!U{b\ll{a}};(4,0.5) *+!D{b}
\ar|(0.5)\hole ,(4,-0.5) *+!U{a\ul{b}};(3,0.5) *+!D{a}
\ar|(0.5)\hole ,(0,-0.5) *+!U{b\ur{a}};(1,0.5) *+!D{b}
\ar ,(1,-0.5)*+!U{a\lr{b}} ;(0,0.5)*+!D{a}
\endxy
\end{equation}
\end{thm}

\begin{proof}
These equations result from the axioms given in~\cite{LouVirt}, which
are required for the biquandle to be invariant under Reidemeister II:
\begin{equation} \label{E:BQRII}
x\ =\ x\lr{y}\Bll{y\ur{x}}
\ =\  x\ur{y}\Bul{y\lr{x}}
\ =\  x\ul{y}\Bur{y\ll{x}}
\ =\  x\ll{y}\Blr{y\ul{x}}
\end{equation}

To derive the ``inverse equations'' for the right handed crossing, start with
the original biquandle equations given in~\eqref{E:BQR}, and solve for
$a$ and $b$ in terms of $X=c$ and $Y=d$:
\begin{align} \label{E:1}
b\ur{a} &= X &\ &\hbox{and} &a\lr{b} &=Y \,.\\
\intertext{Operate on each of these equations appropriately to get}
b\ur{a}\Bul{a\lr{b}} &= X\Bul{a\lr{b}} &\ &\hbox{and} &a\lr{b}\Bll{b\ur{a}} &= Y\Bll{b\ur{a}} ,\notag \\
\intertext{which by~\eqref{E:BQRII}, simplifies to}
b &= X\Bul{a\lr{b}} &\ &\hbox{and} &a &= Y\Bll{b\ur{a}} .\notag \\
\intertext{Finally, backsubstitute from~\eqref{E:1} to get the desired result:}
b &= X\ul{Y} &\ &\hbox{and} &a &= Y\ll{X} .\label{E:BQRrh}
\end{align}
To arrive at the ``inverse equations'' for the left handed
crossing in~\eqref{E:InvABQR}, apply the same argument.
\end{proof}
The five basic morphisms used to compute a knot biquandle are summarized in
Table~\ref{T:BQ_morphisms}.

\begin{table}
    \begin{align*}
        {\xy 0;/r1.25pc/:
            \xoverv~{(-1,1)*+!D{b\ur a}}{(1,1)*+!D{a\lr b}}{(-1,-1)*+!U{a}}{(1,-1)*+!U{b}}=<
        \endxy} &: \ \phi_u(a,b) = (b\ur a, a\lr b) &
        {\xy 0;/r1.25pc/:
            \xunderv~{(-1,1)*+!D{b\ll a}}{(1,1)*+!D{a\ul b}}{(-1,-1)*+!U{a}}{(1,-1)*+!U{b}}=<
        \endxy} &: \ \phi_u^{-1}(a,b) = (b\ll a, a\ul b) \\ \\
        {\xy 0;/r1.25pc/:
            \xunderv~{(-1,1)*+!D{b\lr a}}{(1,1)*+!D{a\ur b}}{(-1,-1)*+!U{a}}{(1,-1)*+!U{b}}=>
        \endxy} &: \ \phi_d^{-1}(a,b) = (b\lr a, a\ur b) &
        {\xy 0;/r1.25pc/:
            \xoverv~{(-1,1)*+!D{b\ul a}}{(1,1)*+!D{a\ll b}}{(-1,-1)*+!U{a}}{(1,-1)*+!U{b}}=>
        \endxy} &: \ \phi_d(a,b) = (b\ul a, a\ll b)
    \end{align*}
    \[
        {\xy 0;/r1.25pc/:0*{\cir<3pt>{}}
            ,\ar@{-} (-1,-1)*+!U{a};(1,1)*+!D{a} *+!D{\phantom{a\lr b}}
            ,\ar@{-} (-1,1)*+!D{b};(1,-1)*+!U{b} *+!U{\phantom{b\ur b}}
        \endxy} : \quad \tau(a,b) = (b,a)
    \]
    \caption{Morphisms used to compute a knot biquandle} \label{T:BQ_morphisms}
\end{table}

One result of the inverse equations is that biquandles will \emph{never} be
able to distinguish the virtual Hopfs (and many other virtual knots)
regardless of how we define the operations.  To see this, view the
downward pointing negative crossing as an automorphism, $\phi_d^{-1}$:
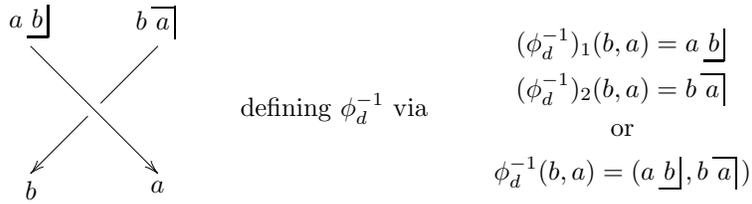
\begin{equation*}
\xy 0;/r2pc/:
    (-1,1)="a" *+!D{a\lr b} ,(1,1)="b" *+!D{b\ur a}
    ,(-1,-1)="c" *+!U{b} ,(1,-1)="d" *+!U{a}
    \xoverv~{"c"}{"d"}{"a"}{"b"}=<
\endxy
\text{\qquad defining $\phi_d^{-1}$ via \qquad}
\begin{gathered}
(\phi_d^{-1})_1(b,a) = a\lr b \\
(\phi_d^{-1})_2(b,a) = b\ur a \\
\text{or}\\
\phi_d^{-1}(b,a) = (a\lr b, b\ur a)
\end{gathered}
\end{equation*}
We've swapped the variables above to demonstrate the similarities between
$\phi_d^{-1}$ and $\phi_u$:
\begin{equation} \label{E:UpVsDown}
\phi_u(a,b)=\phi_d^{-1}(b,a)
\end{equation}
This time, the subscript of $d$ indicates that the direction of the morphism
is downwards, or against the orientation of the strands.  It is indicated as
an inverse morphism because the crossing orientation is negative.

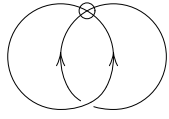
\def \VHopf#1{\xy 0;/r20pt/:0 *\cir(1,0){u^r},(1,0)*\cir(1,0){r^d},a(60)*\cir<3pt>{},#1 \endxy}
\[
\VHopf{/u1pc/,(0,.08)*\dir{>},(1,0.08)*\dir{>}\vunder~{(1,-1)}{(0,-1)}{(1,0)}{(0,0)}}
\]

Let's compute the description for a biquandle associated with $VHopf_{+}$.
The composition $\phi_u\tau$ (see~\ref{T:BQ_morphisms}) represents
the tangle of a positive crossing placed above a virtual crossing.
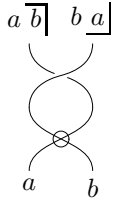
\[
{\xy 0;/r2pc/:
\vtwist~{(0,1)*+!D{a \ur b}}{(1,1)*+!D{b\lr a}}{(0,0)}{(1,0)}
,(0,0);(1,-1) *+!U{b} **\crv{(0,-0.5)&(1,-0.5)},(0,-1)*+!U{a};(1,0) **\crv{(0,-0.5)&(1,-0.5)}
,(0.5,-0.5) *{\cir<3pt>{}}
\endxy}
\]
Connecting the
inputs below to the outputs above gives us a diagram of $VHopf_+$ whose
associated biquandles have the following description:
\[ \bracket{ a,b \ |\ a=a\ur b,\ b=b\lr a }. \]
To compute a biquandle of $VHopf_-$, label the strands across the top with the
generators $a$ on the left and $b$ on the right.
Apply the morphism $\tau\phi_d^{-1}$ on $(b,a)$ to get the description:
\[ \bracket{ b,a \ |\ a=a\ur b,\ b=b\lr a }. \]
Both virtual Hopf links have the same biquandle descriptions.

\begin{figure}
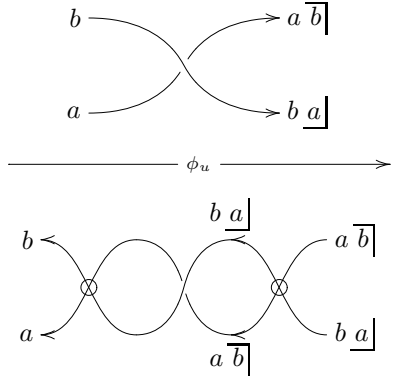

\[{\xy 0;/r3pc/:
 (0,0)="a1"*+!R{a},(2,0)="a2"*+!L{b\lr a}
,(0,1)="b1"*+!R{b},(2,1)="b2"*+!L{a\ur b}
,\htwist~{"b1"}{"b2"}{"a1"}{"a2"}=>
\endxy}\]
\[{\xy 0;/r3pc/:\ar|<(.5){\ \phi_u\ } (0,0);(4,0)\endxy}\]
\[{\xy 0;/r3pc/:
 (0,0)="a1"*+!R{a},(1,0)="a2",(2,0)="a3"*+!U{a\ur b},(3,0)="a4"*+!L{b\lr a}
,(0,1)="b1"*+!R{b},(1,1)="b2",(2,1)="b3"*+!D{b\lr a},(3,1)="b4"*+!L{a\ur b}
,\htwistneg~{"b2"}{"b3"}{"a2"}{"a3"}
,{\save "b1"="temp","a1";"a2":"temp"::(0,0);(1,1) **\crv{(0.5,0)&(0.5,1)}?<(0)*\dir{<},(0,1);(1,0) **\crv{(0.5,1)&(0.5,0)}?<(0)*\dir{<},(0.5,0.5) *{\cir<3pt>{}} \restore}
,{\save "b3"="temp","a3";"a4":"temp"::(0,0);(1,1) **\crv{(0.5,0)&(0.5,1)}?<(0)*\dir{<},(0,1);(1,0) **\crv{(0.5,1)&(0.5,0)}?<(0)*\dir{<},(0.5,0.5) *{\cir<3pt>{}} \restore}
\endxy}\]
\caption{Equivalence Under AD Inversion\label{Fi:AD_Inv_Replacement}}
\end{figure}

\begin{defn} \label{D:ADinv}
If $K$ is an virtual knot diagram, then the \emph{AD inversion of $K$} is obtained by
replacing all classical crossings by their switched virtualizations and reversing the
orientation of the resulting diagram, as shown in Figure~\ref{Fi:AD_Inv_Replacement}.
\end{defn}

\begin{thm} \label{T:BQ_ADinv}
The biquandle is invariant under AD inversion.
\end{thm}
\begin{proof}
Consider equation~\eqref{E:UpVsDown}.  We modify it slightly to get
\begin{subequations} \label{E:ADInv}
\begin{align}
\phi_u&=\tau\phi_d^{-1}\tau \label{E:ADInv1} \\
\phi_u^{-1}&=\tau\phi_d\tau. \label{E:ADInv2}
\end{align}
\end{subequations}
The example crossing replacement in Figure~\ref{Fi:AD_Inv_Replacement} is the personification
of~\eqref{E:ADInv1}. Note that if we perform this replacement on all classical crossings
in a given diagram, then the orientations will all be compatible, and the result will follow.

\end{proof}
\section{Using Braids to Compute the Biquandle} \label{S:BQsViaBraids}
Recall the braid group on $n$ strings.
\begin{thm}[Artin]
The \emph{$n$-string braid group} $\B_n$ has the following presentation:
\begin{equation}
\begin{aligned}
\bigl<
 \sigma_{1},\dots,\sigma_{n-1}\ |
 &\ \sigma_{i}\sigma_{k} = \sigma_{k}\sigma_{i}\ &&|i-k|>1, \\
 &\ \sigma_{i}\sigma_{i+1}\sigma_{i} = \sigma_{i+1}\sigma_{i}\sigma_{i+1}
\bigr>
\end{aligned}
\end{equation}
\end{thm}
For a proof, see~\cite{Birman}.

Each of the generators $\sigma_i$ along with their inverses $\sigma_i^{-1}$ can
be represented by the following diagrams:
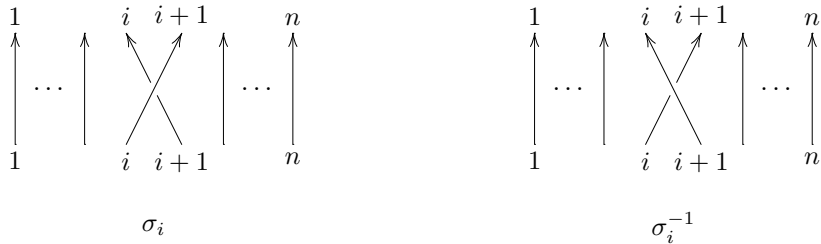
\[ \label{E:BraidGens}
\hbox{

\xy 0;/r1.75pc/:(2.5,-2.5)*{\sigma_i}, (0.6,0) *{\dots}, (4.35,0) *{\dots}, 0
    \xunderv~{(2,-1)*+!U{i}}{(3,-1)*+!U{i+1}}{(2,1)*+!D{i}}{(3,1)*+!D{i+1}}=>
    ,\ar (0,-1)*+!U{1}; (0,1)*+!D{1}
    ,\ar (1.25,-1); (1.25,1)
    ,\ar (3.75,-1); (3.75,1)
    ,\ar (5,-1)*+!U{n}; (5,1)*+!D{n}
\endxy
\qquad \qquad \qquad \qquad
\xy 0;/r1.75pc/:(2.5,-2.5)*{\sigma_i^{-1}}, (0.6,0) *{\dots}, (4.35,0) *{\dots}, 0
    \xoverv~{(2,-1)*+!U{i}}{(3,-1)*+!U{i+1}}{(2,1)*+!D{i}}{(3,1)*+!D{i+1}}=>
    ,\ar (0,-1)*+!U{1}; (0,1)*+!D{1}
    ,\ar (1.25,-1); (1.25,1)
    ,\ar (3.75,-1); (3.75,1)
    ,\ar (5,-1)*+!U{n}; (5,1)*+!D{n}
\endxy
}
\]
Add a further generator $\nu_i$, corresponding to a virtual crossing
between the $i^{th}$ and $(i+1)^{th}$ strands to generalize the
$n$-strand Braid group to virtual braids.
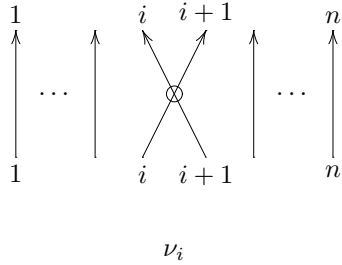
\[
\xy 0;/r2pc/:(2.5,-2.5)*{\nu_i}, (2.5,0)*\cir<3pt>{}, (0.6,0)*{\dots}, (4.35,0)*{\dots}
    ,\ar (2,-1)*+!U{i}; (3,1)*+!D{i+1}
    ,\ar (3,-1)*+!U{i+1}; (2,1)*+!D{i}
    ,\ar (0,-1)*+!U{1}; (0,1)*+!D{1}
    ,\ar (1.25,-1); (1.25,1)
    ,\ar (3.75,-1); (3.75,1)
    ,\ar (5,-1)*+!U{n}; (5,1)*+!D{n}
\endxy
\]
In this way, the $n$ string braid group extends to virtual theory
(see~\cite{Kamada, KamadaVirtualBraids, LouSofia, LouHellas}).

\begin{defn}
The \emph{$n$-string virtual braid group} $\V\B_n$ has the following presentation:
\begin{equation*}
\begin{aligned}
\bigl<
 \sigma_{1},\dots,\sigma_{n-1} ,\nu_{1},\dots,\nu_{n-1}\ |
 &\ \sigma_{i}\sigma_{k} = \sigma_{k}\sigma_{i} && |i-k|>1, \\
 &\ \sigma_{i}\sigma_{i+1}\sigma_{i} = \sigma_{i+1}\sigma_{i}\sigma_{i+1} \\ 
 &\ \nu_i^2=1, \\ 
 &\ \nu_{i}\nu_{k} = \nu_{k}\nu_{i} && |i-k|>1, \\
 &\ \nu_{i}\nu_{i+1}\nu_{i} = \nu_{i+1}\nu_{i}\nu_{i+1} \\ 
 &\ \sigma_{i}\nu_{k} = \nu_{k}\sigma_{i} && |i-k|>1, \\
 &\ \sigma_{i}\nu_{i+1}\nu_{i} = \nu_{i+1}\nu_{i}\sigma_{i+1} 
\bigr>
\end{aligned}
\end{equation*}
\end{defn}
In~\cite{Fenn, FennEtal}, Fenn, Rim\`anyi, and Rourke prove that the welded
braid group on $n$ strands is isomorphic to the automorphism group of the
free quandle\footnote{Recall that welded braids are virtual braids with the
addition of an extra mixed move.} on $n$ generators. The morphisms described
above make it possible to analyze the automorphism group of a biquandle.
A virtual braid determines a morphism of biquandle elements.  Their result
motivated the approach in this paper.

\begin{defn}
For any braid $\beta\in\V\B_n$, the closure $cl(\beta)$ is a knot obtained by identifying
initial points and end points of each of the braid strings
\end{defn}

We will need a generalization of Alexander's Theorem to virtual braids,
proved independently by the authors and by Kamada~\cite{Kamada}.

\begin{thm}[Alexander's Theorem for Virtual Braids] \label{L:VirtualAlexander}
If $K$ is a virtual knot, then for some $n$, there exists a braid $\beta\in\V\B_n$
such that the closure of $\beta$ is $K$.
\end{thm}

\begin{proof}
Suppose $K$ is a virtual knot.  Choose a representative knot diagram and set a
braid axis in a arbitrary region in the diagram.  Now, choose an orientation for
the knot, and rotate each real crossing via planar isotopy so that both the over
and under crossings align in a clockwise direction around the axis.  Note that
this leaves a finite number of virtual arcs connecting each real crossing to another,
and it is only these connecting arcs which might not travel counter-clockwise
around the axis.  One by one, we can perform a detour move on each of these arcs
so that they travel clockwise between each real crossing.  After performing all of
these detour moves, we will have a virtual knot diagram for which all the real
crossings are aligned clockwise around the braid axis and every connecting path
between these crossings is also aligned clockwise around the axis.  Since this
new diagram has only a finite number of real and virtual crossings, there exists
a ray from the axis to the exterior of the diagram which avoids every crossing.
Cut along this ray to get the braid $B$ whose closure is equivalent to $K$.
\end{proof}

\begin{thm} \label{T:BQ_BraidInv}
If $BQ$ is a biquandle, then for every $\beta\in \V\B_n$, $BQ(cl(\beta))\cong BQ(cl(\beta^{-1}))$.
\end{thm}

\begin{proof}
Let $\beta = x_{1}x_{2}\cdots x_{k}\in \V\B_n$,
where
\[
x_{i}\in \{ \sigma_{1}^{\epsilon_1},\dots,\sigma_{(n-1)}^{\epsilon_{(n-1)}},
 \nu_{1},\dots,\nu_{(n-1)}\ |\ \epsilon_i=\pm 1\}\subset \V\B_n.
\]
Define a map $f_u^{(n)}:\V\B_n \to Aut(BQ^n)$ by setting
\begin{alignat*}{2}
f_u^{(n)}(\sigma_i^{\epsilon}) &= \phi_{u,i}^{(n)} & &= 1_{i-1}\oplus\phi_u^{\epsilon}\oplus 1_{n-1-i} \\
f_u^{(n)}(\nu_i^{\epsilon})    &= \tau_{i}^{(n)} & &= 1_{i-1}\oplus\tau\oplus 1_{n-1-i} \\
\end{alignat*}
and extending $f_u^{(n)}$ to all of $\V\B_n$ using the multiplication reversing rule:
\begin{equation}
f_u^{(n)}(xy) = f_u^{(n)}(y)f_u^{(n)}(x).
\end{equation}
We will use $f_u^{(n)}$ to construct a description for $BQ(cl(\beta))$.
Label the bottom strands of $\beta$ with the generators $a_1, a_2, \dots, a_n$ ordered
from left to right.  To get the equations, apply the composition of
upward morphisms described by $\beta$ and equate the resulting labels on the
upper strands with the corresponding generators below.  In terms of $f_u^{(n)}$,
the description is
\begin{equation} \label{E:P4beta}
BQ(cl(\beta)) = \bracket{a_1, a_2, \dots, a_n\ |\ f_u^{(n)}(a_1, a_2, \dots, a_n)=(a_1, a_2, \dots, a_n)}
\end{equation}

Define the map $f_d^{(n)}:\V\B_n \to Aut(BQ^n)$ by setting
\begin{alignat*}{2}
f_d^{(n)}(\sigma_i^{\epsilon}) &= \phi_{u,i}^{(n)} & &= 1_{n-1-i}\oplus\phi_d^{-\epsilon}\oplus 1_{i-1} \\
f_d^{(n)}(\nu_i^{\epsilon})    &= \tau_{i}^{(n)}     & &= 1_{n-1-i}\oplus\tau\oplus 1_{i-1} \\
\end{alignat*}
and extending $f_d^{(n)}$ to all of $\V\B_n$ using the multiplication preserving rule:
\begin{equation}
f_d^{(n)}(xy) = f_d^{(n)}(x)f_d^{(n)}(y).
\end{equation}
We use $f_d^{(n)}$ to construct a description for $BQ(cl(\beta^{-1}))$.
Label the top strands of $\beta^{-1}$ with the generators $a_1, a_2, \dots, a_n$
ordered from left to right.  Recall this is the reverse of the counter-clockwise
convention for inputs.  The description is:
\begin{equation} \label{E:P4betaInv}
BQ(cl(\beta^{-1})) = \bracket{a_1, a_2, \dots, a_n\ |\ f_d^{(n)}(a_n, \dots, a_2, a_1)=(a_n, \dots, a_2, a_1)}
\end{equation}

Express $f_u^{(n)}(\beta)$ and $f_d^{(n)}(\beta^{-1})$ in terms of the $x_i$'s:
\begin{align*}
f_u^{(n)}(\beta)  &= f_u^{(n)}(x_{1}x_{2}\cdots x_{k}) \\
            &= f_u^{(n)}(x_{k})\cdots f_u^{(n)}(x_{2})f_u^{(n)}(x_{1}) \\
\\
f_d^{(n)}(\beta^{-1}) &= f_d^{(n)}(x_{k}^{-1}\cdots x_{2}^{-1}x_{1}^{-1}) \\
                &= f_d^{(n)}(x_{k}^{-1})\cdots f_d^{(n)}(x_{2}^{-1})f_d^{(n)}(x_{1}^{-1}) \\
\end{align*}
We claim that these compositions generate the same set of equations
for the descriptions given in~\eqref{E:P4beta} and in~\eqref{E:P4betaInv}.
Start with the virtual morphisms:
\begin{align*}
f_u^{(n)}(\nu_i)(a_1, \dots, a_i, a_{(i+1)}, \dots, a_n)
  &= (a_1, \dots, a_{(i+1)}, a_i, \dots, a_n) \\
f_d^{(n)}(\nu_i^{-1})(a_n, \dots, a_{(i+1)},  a_i, \dots, a_1)
  &= (a_n, \dots, a_i, a_{(i+1)}, \dots, a_1).
\end{align*}
Both of the above permute $a_i$ with $a_{(i+1)}$.  Moreover they are reversals of each other,
\begin{equation} \label{E:T_u_vs_d}
f_u^{(n)}(\nu_i) = T_n f_d^{(n)}(\nu_i^{-1})
\end{equation}
where $T_n$ is the map which reverses $n$-tuples:
$T_n(b_1, b_2, \dots, b_n)=(b_n, \dots, b_2, b_1)$.

Recall the relationship between $\phi_{u}$ and $\phi_{d}^{-1}$.
Equation~\eqref{E:UpVsDown} can be generalized for $\phi_{(u,i)}^{\epsilon}$
and $\phi_{(d,i)}^{-\epsilon}$.
Consider how it works for a specific example: the braids $\sigma_2$
and $\sigma_2^{-1}$ in $\V\B_3$.
The morphism $f_u^{(3)}(\sigma_2)=\phi_{(u,i)}$ mapping upwards on the page looks like

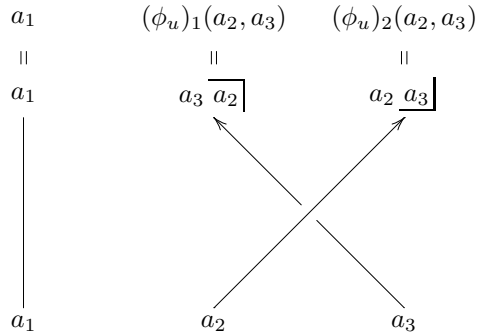
\[
\xy 0;/r6pc/:0*+!U{a_1};(0,1) **\dir{-}
,(0,1)="t"; "t"+/u1pt/, "t" + (0,0.12) *{a_1} + (0,0.16) *\dir{=} + (0,0.24) *{a_1}
,(1,1)="t"; "t"+/u1pt/, "t" + (0,0.12) *{a_3 \ur{a_2}} + (0,0.16) *\dir{=} + (0,0.24) *{(\phi_u)_1(a_2,a_3)}
,(2,1)="t"; "t"+/u1pt/, "t" + (0,0.12) *{a_2 \lr{a_3}} + (0,0.16) *\dir{=} + (0,0.24) *{(\phi_u)_2(a_2,a_3)}
\ar (1,0)*+!U{a_2}; (2,1)
,\ar|(0.5)\hole (2,0)*+!U{a_3}; (1,1)
\endxy
\]
and the morphism $f_d^{(3)}(\sigma_2^{-1})=\phi_{(d,i)}^{-1}$ mapping downwards on the page looks like

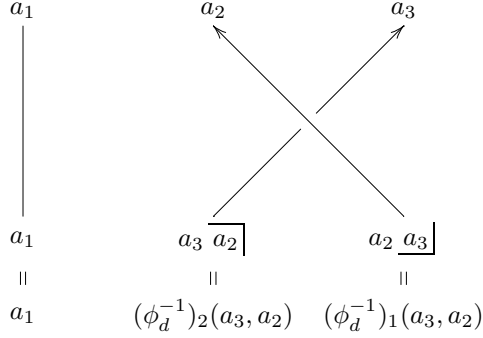
\[
\xy 0;/r6pc/:0;(0,1)*+!D{a_1} **\dir{-}
,(0,0)="t"; "t"-/u1pt/, "t" - (0,0.12) *{a_1} - (0,0.24) *\dir{=} - (0,0.16) *{a_1}
,(1,0)="t"; "t"-/u1pt/, "t" - (0,0.12) *{a_3 \ur{a_2}} - (0,0.24) *\dir{=} - (0,0.16) *{(\phi_d^{-1})_2(a_3,a_2)}
,(2,0)="t"; "t"-/u1pt/, "t" - (0,0.12) *{a_2 \lr{a_3}} - (0,0.24) *\dir{=} - (0,0.16) *{(\phi_d^{-1})_1(a_3,a_2)}
\ar|(0.5)\hole (1,0); (2,1)*+!D{a_3}
,\ar (2,0); (1,1)*+!D{a_2}
\endxy
\]

The outputs of both maps read the same from left to right.
This means that the $3$-tuples $\phi_{u,2}(a_1,a_2,a_3)$ and
$\phi_{d,2}^{-1}(a_3,a_2,a_1)$ are reversals of each other.
More generally,
\begin{equation} \label{E:phi_u_vs_d}
   f_u^{(n)}(\sigma_i^\epsilon)(a_1, a_2, \dots, a_n)
   = T_n f_d^{(n)}(\sigma_i^{-\epsilon})(a_n, \dots, a_2, a_1).
\end{equation}

Equations~\eqref{E:T_u_vs_d} and~\eqref{E:phi_u_vs_d} imply that
$f_u^{(n)}(x_{i})=T_n f_d^{(n)}(x_{i}^{-1})$ for each $i$.
The result carries through to compositions.  If we equate
\begin{align*}
(a_1,a_2,\dots ,a_n)
&=f_u^{(n)}(x_{k})\cdots f_u^{(n)}(x_{2})f_u^{(n)}(x_{1})(a_1,a_2,\dots ,a_n) \\
\intertext{then it must also be true that}
T_n(a_n,\dots a_2,a_1)
&= T_n f_d^{(n)}(x_{k}^{-1})\cdots f_d^{(n)}(x_{2}^{-1})f_d^{(n)}(x_{1}^{-1})(a_n,\dots ,a_2,a_1).
\end{align*}
Therefore, the descriptions in~\eqref{E:P4beta} and~\eqref{E:P4betaInv} are the same.
This concludes our proof of the theorem.
\end{proof}

\noindent{\bf Remark.}
We can summarize this argument as follows.  Consider two braids, $\beta$ and $\beta^{-1}$:

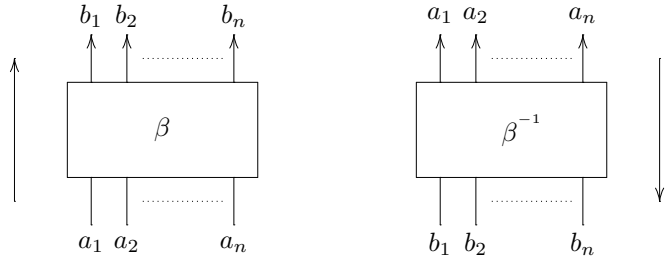
\def \BraidBox#1#2#3{\xy 0*\txt{#1};/r1.5pc/:
(-2,-1)="a",(-2,1)="b",(2,1)="c",(2,-1)="d",(0,.4)="dy",
@={"a","b","c","d"},
s0="prev" @@{;"prev";**\dir{-}="prev"},
{\ar (-1.50,1);(-1.50,2)="u1"},
{\ar (-0.75,1);(-0.75,2)="u2"},
{\ar (1.50,1);(1.50,2)="u3"},
{\ar@{.} (-0.5,1.5);(1.25,1.5)},
{\ar@{-} (-1.50,-1);(-1.50,-2)="d1"},
{\ar@{-} (-0.75,-1);(-0.75,-2)="d2"},
{\ar@{-} (1.50,-1);(1.50,-2)="d3"},
{\ar@{.} (-0.5,-1.5);(1.25,-1.5)},
 "u1"+"dy" *+{#2_1},
 "u2"+"dy" *+{#2_2},
 "u3"+"dy" *+{#2_n},
 "d1"-"dy" *+{#3_1},
 "d2"-"dy" *+{#3_2},
 "d3"-"dy" *+{#3_n}
\endxy}

\[
\xy 0;/r1.5pc/:0; \ar (0.0,-1.5);(0.0,1.5)\endxy
\qquad
\BraidBox{$\beta$}{b}{a}
\qquad\qquad\qquad
\BraidBox{$\!\quad\beta^{^{-1}}$}{a}{b}
\qquad
\xy 0;/r1.5pc/:0; \ar (0.0,1.5);(0.0,-1.5)\endxy
\]
If we view $\beta$ as an upward
morphism, and $\beta^{-1}$ as a downward morphism, then it is clear that the two morphisms
give rise to identical biquandle descriptions.
The crucial step in seeing this is illustrated in the final two diagrams in the above proof.

Recall that the fundamental group $\pi_1(K)$ is a quandle.
One consequence of Theorem~\ref{T:BQ_BraidInv} is the following well known
result, generalized to virtual braids.
\begin{cor} \label{C:BQ_BraidInv1}
If $\beta \in \V\B_n$, $\pi_1(cl(\beta))=\pi_1(cl(\beta^{-1}))$
\end{cor}


\begin{cor} \label{C:BQ_BraidInv2}
If $K$ is a virtual knot diagram and $K^\uparrow$ is the vertical mirror image of $K$,
then $BQ(K)\cong BQ(K^\uparrow)$.
\end{cor}
\begin{proof}
By following the proof of Theorem~\ref{L:VirtualAlexander}, we can simultaneously convert
a knot and its vertical mirror image to closed braid forms that are vertical mirror images
of one another.
\end{proof}

\begin{figure}
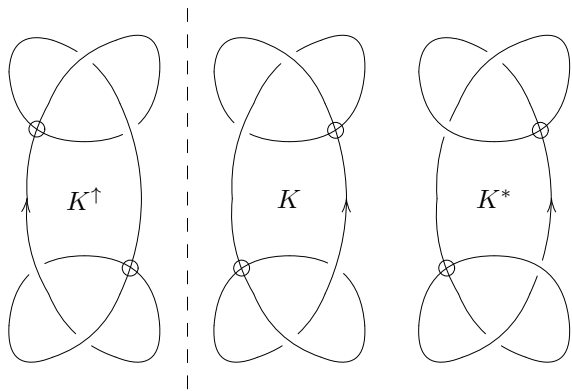

\[
{\xy 0;/u1pc/:
 (0,0)*+{K^{\uparrow}},
 (0.,1.8)="a1", (0,-1.8)="a2",
 (1.8,0.0)="b1",
 (2.,3.0)="c0", (4.,2.375)="c1", (3.,1.125)="c2", (3.,-1.125)="c3", (4.,-2.375)="c4", (2.,-3.0)="c5",
 (5.75,3)="d0", (5.,1.)="d1", (5.,-1.)="d2", (5.75,-3)="d3"
 ,{\ar@{-}@/^0.23pc/ "a1";"c2"}?(.72)*{\cir<3pt>{}}
 ,{\ar@{-}@/^9.pt/ "b1";"c1" }
 ,{\ar@{-}@/^1.5ex/ "c1";"d1"}
 ,{\ar@{-}@/_2.pt/ "a2";"c3"}
 ,{\ar@{-}|<(.405){\hole}@/_9.pt/ "b1";"c4" }
 ,{\ar@{-}@/^8.pt/ "d2";"c4"}
 ,{\ar@{-}|<(.62){\hole}@/_4.pt/ "c3";"d1"}
 ,{\ar@{-}@/^4.pt/ "c2";"d2"}
,
 (0.,1.8)="a1", (0,-1.8)="a2",
 (-1.8,0.0)="b1",
 (-2.,3.0)="c0", (-4.,2.375)="c1", (-3.,1.125)="c2", (-3.,-1.125)="c3", (-4.,-2.375)="c4", (-2.,-3.0)="c5",
 (-5.75,3)="d0", (-5.,1.)="d1", (-5.,-1.)="d2", (-5.75,-3)="d3"
 ,{\ar@{<-}@/_3.pt/ "a1";"c2"}
 ,{\ar@{-}|<(.42){\hole}@/_9.pt/ "b1";"c1" }
 ,{\ar@{-}@/_8.pt/ "c1";"d1"}
 ,{\ar@{-}@/^2.pt/ "a2";"c3"}?(.72)*{\cir<3pt>{}}
 ,{\ar@{-}@/^9.pt/ "b1";"c4" }
 ,{\ar@{-}@/_8.pt/ "d2";"c4"}
 ,{\ar@{-}@/^4.pt/ "c3";"d1"}
 ,{\ar@{-}|<(.62){\hole}@/_4.pt/ "c2";"d2"}
\endxy}
\quad {\xy 0;/r1pc/: (0,-6)="d",(0,6)="u" ,{\ar@{--} "d";"u"}
\endxy}
\quad {\xy 0;/u1pc/:
 (0,0)*+{K},
 (0.,1.8)="a1", (0,-1.8)="a2",
 (1.8,0.0)="b1",
 (2.,3.0)="c0", (4.,2.375)="c1", (3.,1.125)="c2", (3.,-1.125)="c3", (4.,-2.375)="c4", (2.,-3.0)="c5",
 (5.75,3)="d0", (5.,1.)="d1", (5.,-1.)="d2", (5.75,-3)="d3"
 ,{\ar@{-}@/^0.23pc/ "a1";"c2"}
 ,{\ar@{-}|<(.42){\hole}@/^9.pt/ "b1";"c1" }
 ,{\ar@{-}@/^1.5ex/ "c1";"d1"}
 ,{\ar@{-}@/_2.pt/ "a2";"c3"}
 ,{\ar@{-}@/_9.pt/ "b1";"c4" }?(.405)*{\cir<3pt>{}}
 ,{\ar@{-}@/^8.pt/ "d2";"c4"}
 ,{\ar@{-}@/_4.pt/ "c3";"d1"}
 ,{\ar@{-}|<(.62){\hole}@/^4.pt/ "c2";"d2"}
,
 (0.,1.8)="a1", (0,-1.8)="a2",
 (-1.8,0.0)="b1",
 (-2.,3.0)="c0", (-4.,2.375)="c1", (-3.,1.125)="c2", (-3.,-1.125)="c3", (-4.,-2.375)="c4", (-2.,-3.0)="c5",
 (-5.75,3)="d0", (-5.,1.)="d1", (-5.,-1.)="d2", (-5.75,-3)="d3"
 ,{\ar@{-}@/_3.pt/ "a1";"c2"}
 ,{\ar@{-}@/_9.pt/ "b1";"c1" }?(.42)*{\cir<3pt>{}}
 ,{\ar@{-}@/_8.pt/ "c1";"d1"}
 ,{\ar@{<-}@/^2.pt/ "a2";"c3"}
 ,{\ar@{-}|<(.405){\hole}@/^9.pt/ "b1";"c4" }
 ,{\ar@{-}@/_8.pt/ "d2";"c4"}
 ,{\ar@{-}|<(.62){\hole}@/^4.pt/ "c3";"d1"}
 ,{\ar@{-}@/_4.pt/ "c2";"d2"}
\endxy}
\qquad {\xy 0;/u1pc/:
 (0,0)*+{K^{*}},
 (0.,1.8)="a1", (0,-1.8)="a2",
 (1.8,0.0)="b1",
 (2.,3.0)="c0", (4.,2.375)="c1", (3.,1.125)="c2", (3.,-1.125)="c3", (4.,-2.375)="c4", (2.,-3.0)="c5",
 (5.75,3)="d0", (5.,1.)="d1", (5.,-1.)="d2", (5.75,-3)="d3"
 ,{\ar@{-}|<(.72){\hole}@/^0.23pc/ "a1";"c2"}
 ,{\ar@{-}@/^9.pt/ "b1";"c1" }
 ,{\ar@{-}@/^1.5ex/ "c1";"d1"}
 ,{\ar@{>-}@/_2.pt/ "a2";"c3"}
 ,{\ar@{-}@/_9.pt/ "b1";"c4" }?(.405)*{\cir<3pt>{}}
 ,{\ar@{-}@/^8.pt/ "d2";"c4"}
 ,{\ar@{-}|<(.62){\hole}@/_4.pt/ "c3";"d1"}
 ,{\ar@{-}@/^4.pt/ "c2";"d2"}
,
 (0.,1.8)="a1", (0,-1.8)="a2",
 (-1.8,0.0)="b1",
 (-2.,3.0)="c0", (-4.,2.375)="c1", (-3.,1.125)="c2", (-3.,-1.125)="c3", (-4.,-2.375)="c4", (-2.,-3.0)="c5",
 (-5.75,3)="d0", (-5.,1.)="d1", (-5.,-1.)="d2", (-5.75,-3)="d3"
 ,{\ar@{-}@/_3.pt/ "a1";"c2"}
 ,{\ar@{-}@/_9.pt/ "b1";"c1" }?(.42)*{\cir<3pt>{}}
 ,{\ar@{-}@/_8.pt/ "c1";"d1"}
 ,{\ar@{-}|<(.72){\hole}@/^2.pt/ "a2";"c3"}
 ,{\ar@{-}@/^9.pt/ "b1";"c4" }
 ,{\ar@{-}@/_8.pt/ "d2";"c4"}
 ,{\ar@{-}@/^4.pt/ "c3";"d1"}
 ,{\ar@{-}|<(.62){\hole}@/_4.pt/ "c2";"d2"}
\endxy}
\]
\caption{The two mirror images\label{Fi:Mirror}}
\end{figure}

\noindent{\bf Remark.}
Corollary~\ref{C:BQ_BraidInv2} shows that the biquandle is the same for a virtual link and its
vertical mirror image.  In Figure~\ref{Fi:Mirror}, we give an example of virtual knot $K$
and two mirror images: the vertical mirror image $K^\uparrow$, and the mirror image $K^*$,
obtained by switching all the crossings in the plane.
By the Corollary, the biquandle of $K^\uparrow$ is equal to the biquandle of $K$. However, the
biquandle of $K^*$ is not isomorphic to the biquandle of $K$.  In particular, calculation shows
that $\pi_1(K)$ is the same as the fundamental group of the classical trefoil knot, and hence
the biquandle is non-trivial.
Calculation also shows that the fundamental group of $K^*$ is trivial.  Hence, $K$ and $K^*$
have non-isomorphic biquandles.  Furthermore, one sees that $K$ and $K^*$ have the same Jones
polynomial, whereas the Jones polynomial of $K^\uparrow$ is different.

\section{The Alexander Biquandle} \label{S:ABQ}
Consider any module $M$ over the ring $R=\Z[s,s^{-1},t,t^{-1}]$.
Defining the binary operations with the following equations provides us
with a biquandle structure on $M$:
\begin{alignat}{2} \label{E:ABQoperations}
    a\ur b &= ta+(1-st)b    & \qquad  a\lr b &= sa \\
    a\ul b &= \frac{1}{t}a+(1-\frac{1}{st})b    & \qquad   a\ll b &= \frac{1}{s}a
\end{alignat}
If $M$ is a free module, we call this a \emph{free Alexander biquandle}.

We associate a specific biquandle to a virtual knot diagram by
taking the free Alexander module obtained by assigning one generator
for each arc and factoring out by the submodule generated by the
relations given in~\eqref{E:ABQoperations}. We call the resulting
biquandle $ABQ(K)$ the \emph{Alexander biquandle of the virtual knot~$K$}.

\begin{figure}
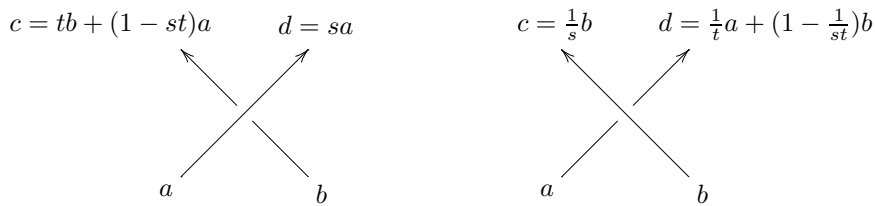

\[
\xy 0;<4pc,0pc>:
 (0,0.5)="ul1", (1,0.5)="ur1", (3,0.5)="ul2", (4,0.5)="ur2",
(0,-0.5)="dl1",(1,-0.5)="dr1",(3,-0.5)="dl2",(4,-0.5)="dr2",
(0.3,0)="dx", (0,0.2)="dy",
"dl1" *+!UR{a}, "dr1" *+!UL{b},
"ul1"+"dx"+"dy" *+!R{c=tb+(1-st)a}, "ur1"-"dx"+"dy" *+!L{d=sa},
"dl2" *+!UR{a}, "dr2" *+!UL{b},
"ul2"+"dx"+"dy" *+!R{c=\frac{1}{s}b}, "ur2"-"dx"+"dy" *+!L{d=\frac{1}{t}a+(1-\frac{1}{st})b},
\ar , "dl1"; "ur1"
\ar|(0.5)\hole , "dr1"; "ul1"
\ar|(0.5)\hole , "dl2"; "ur2"
\ar , "dr2"; "ul2"
\endxy
\]
\caption{\label{Fi:ABQrelations} The Alexander biquandle relations}
\end{figure}
The biquandle morphisms are linear maps, listed in~\ref{T:ABQMatrices}.
In the context of the Alexander biquandle, the morphisms from before become
\[
        \phi_u=A \qquad \phi_d=\widehat{B}
 \qquad \phi_u^{-1}=B \qquad \phi_d^{-1}=\widehat{A}
 \qquad \tau=V
\]
We have included four more matrices in Table~\ref{T:ABQMatrices} in addition to the morphisms
given in Table~\ref{T:BQ_morphisms}.  They are
the linear maps which result when a change of basis transforms the upward and downward
morphisms into left and right morphisms.
\begin{table}
    \begin{align*}
        \xy 0;/r3pc/:(0,0.5) \xoverv=< \endxy: \quad A &=
            \begin{pmatrix} 1-st&t \\ s&0 \end{pmatrix} &
        \xy 0;/r3pc/:(0,0.5) \xunderv=> \endxy : \quad \widehat{A} &=
            \begin{pmatrix} 0&s \\ t&1-st \end{pmatrix} \\ \\
        \xy 0;/r3pc/:(0,0.5) \xunderv=< \endxy : \quad B &=
            \begin{pmatrix} 0&\frac{1}{s} \\ \frac{1}{t}&1-\frac{1}{st} \end{pmatrix} &
        \xy 0;/r3pc/:(0,0.5) \xoverv=> \endxy : \quad \widehat{B} &=
            \begin{pmatrix} 1-\frac{1}{st}&\frac{1}{t} \\ \frac{1}{s}&0\end{pmatrix} \\ \\
        \xy 0;/r3pc/:(0,0.5) \xoverh=> \endxy : \quad C &=
            \begin{pmatrix} 0&\frac{1}{s} \\ t&\frac{1}{s}-t \end{pmatrix} &
        \xy 0;/r3pc/:(0,0.5) \xunderh=> \endxy : \quad \widehat{C} &=
            \begin{pmatrix} \frac{1}{s}-t&t \\ \frac{1}{s}&0 \end{pmatrix} \\ \\
        \xy 0;/r3pc/:(0,0.5) \xoverh=< \endxy : \quad D &=
            \begin{pmatrix} 0&s \\ \frac{1}{t}&s-\frac{1}{t} \end{pmatrix} &
        \xy 0;/r3pc/:(0,0.5) \xunderh=< \endxy : \quad \widehat{D} &=
            \begin{pmatrix} s-\frac{1}{t}&\frac{1}{t} \\ s&0 \end{pmatrix}
    \end{align*}
    \bigskip
    \[
        \xy 0;/r1.5pc/:0*{\cir<3pt>{}},\ar@{-} (-1,-1);(1,1),\ar@{-} (-1,1);(1,-1)
        \endxy : \quad V = \begin{pmatrix} 0&1 \\ 1&0 \end{pmatrix}
    \]
    \caption{Matrices for the Alexander biquandle}
    \label{T:ABQMatrices}
\end{table}

There are advantages to considering the Alexander biquandle.
The binary operations are fairly simple, which allows more flexibility
when choosing the morphisms.
In addition, the set of relations for a presentation of $ABQ(K)$
contains a generalization the Alexander polynomial
(see~\cite{Jaeger,LouRad,SawollekACPolys,SilverWilliams}).
\begin{defn}
The \emph{Generalized Alexander Polynomial of $K$,} $G_K(s,t)$ is the determinant
of the relation matrix from a presentation of $ABQ(K)$.  Up to multiples of $\pm s^i t^j$
for $i,j\in\Z$, it is an invariant of $K$.
\end{defn}

A significant problem to study comes from our general results on
biquandle equivalent pairs of knots, or rather pairs of knots
such as those arising from Theorem~\ref{T:BQ_ADinv} and
Theorem~\ref{T:BQ_BraidInv}.
We show two such pairs in Figure~\ref{Fi:VKPairs}.  These
pairs both come from an application of Theorem~\ref{T:BQ_ADinv}.
Relation matrices for the Alexander modules of each pair are also
given in the Figure,
and we remark that the generalized Alexander polynomials are
the same for all four knots, not just between the pairs: $(s-1)(t-1)(st-1)$.

It is unclear why $G_K(s,t)$ cannot distinguish between the two pairs.
In fact, by adding more twists to the vertical tangle, we can extend
the two virtual knot pairs to an infinite set of virtual knot pairs
with the same generalized Alexander polynomials.
Note that it still is possible for the Alexander modules of the
resulting pairs to be different.  We conjecture that this is the case.

\begin{figure}
$$\xy 0;/r1pc/:(0,2)="a1",(1,2)="b1",(0,1)="a2",(1,1)="b2"
    ,(0,0)="a3",(1,0)="b3",(0,-1)="a4",(1,-1)="b4"
    ,(-.5,-1.5)="a5",(1.5,-1.5)="b5",(0.5,-1.5)="c5",(0,-2)="a6",(1,-2)="b6"
    ,(0.5,1.80)*{\cir<3pt>{}},\ar@/^3pt/@{-} "a1";"b2",\ar@/_3pt/@{-} "b1";"a2"
    \vtwist~{"a2"}{"b2"}{"a3"}{"b3"}\vtwist~{"a3"}{"b3"}{"a4"}{"b4"}
    \xunderv~{"a5"}{"a4"}{"a6"}{"c5"}\xunderv~{"b6"}{"c5"}{"b5"}{"b4"}
    \hcap~{"b1"}{"b1"+(3,0.05)}{"b5"}{"b5"+(2,-.05)}
    \hcap~{"a1"}{"a1"-(3,-0.05)}{"a5"}{"a5"-(2,0.05)}|>
    \vcap~{"a6"-(0,1)}{"b6"-(0,1)}{"a6"}{"b6"}
\endxy
\xy 0;/r1pc/:(0,2)="a1",(1,2)="b1",(0,1)="a2",(1,1)="b2"
    ,(0,0)="a3",(1,0)="b3",(0,-1)="a4",(1,-1)="b4"
    ,(-.5,-1.5)="a5",(1.5,-1.5)="b5",(0.5,-1.5)="c5",(0,-2)="a6",(1,-2)="b6"
    ,(0.5,1.80)*{\cir<3pt>{}},\ar@/^3pt/@{-} "a1";"b2",\ar@/_3pt/@{-} "b1";"a2"
    \vtwistneg~{"a2"}{"b2"}{"a3"}{"b3"}\vtwistneg~{"a3"}{"b3"}{"a4"}{"b4"}
    \xoverv~{"a5"}{"a4"}{"a6"}{"c5"}\xoverv~{"b6"}{"c5"}{"b5"}{"b4"}
    \hcap~{"b1"}{"b1"+(3,0.05)}{"b5"}{"b5"+(2,-.05)}
    \hcap~{"a1"}{"a1"-(3,-0.05)}{"a5"}{"a5"-(2,0.05)}|>
    \vcap~{"a6"-(0,1)}{"b6"-(0,1)}{"a6"}{"b6"}
\endxy
\quad
\begin{vmatrix}
    -1&0&-t(st-1)&s^2t^2-st+1\\
    0&-1&st&-s(st-1)\\
    -1&st-1&-st+2&0\\
    0&st&-st+1&-1
\end{vmatrix}
$$
\bigskip
$$
\xy 0;/r1pc/:(0,3)="a1",(1,3)="b1",(0,2)="a2",(1,2)="b2"
    ,(0,1)="a3",(1,1)="b3",(0,0)="a4",(1,0)="b4"
    ,(0,-1)="a5",(1,-1)="b5",(0,-2)="a6",(1,-2)="b6"
    ,(-.5,-2.5)="a7",(1.5,-2.5)="b7",(0.5,-2.5)="c7",(0,-3)="a8",(1,-3)="b8"
    ,(0.5,2.80)*{\cir<3pt>{}},\ar@/^3pt/@{-} "a1";"b2",\ar@/_3pt/@{-} "b1";"a2"
    \vtwist~{"a2"}{"b2"}{"a3"}{"b3"}\vtwist~{"a3"}{"b3"}{"a4"}{"b4"}
    \vtwist~{"a4"}{"b4"}{"a5"}{"b5"}\vtwist~{"a5"}{"b5"}{"a6"}{"b6"}
    \xunderv~{"a7"}{"a6"}{"a8"}{"c7"}\xunderv~{"b8"}{"c7"}{"b7"}{"b6"}
    \hcap~{"b1"}{"b1"+(3,0.05)}{"b7"}{"b7"+(2,-.05)}
    \hcap~{"a1"}{"a1"-(3,-0.05)}{"a7"}{"a7"-(2,0.05)}|>
    \vcap~{"a8"-(0,1)}{"b8"-(0,1)}{"a8"}{"b8"}
\endxy
\xy 0;/r1pc/:(0,3)="a1",(1,3)="b1",(0,2)="a2",(1,2)="b2"
    ,(0,1)="a3",(1,1)="b3",(0,0)="a4",(1,0)="b4"
    ,(0,-1)="a5",(1,-1)="b5",(0,-2)="a6",(1,-2)="b6"
    ,(-.5,-2.5)="a7",(1.5,-2.5)="b7",(0.5,-2.5)="c7",(0,-3)="a8",(1,-3)="b8"
    ,(0.5,2.80)*{\cir<3pt>{}},\ar@/^3pt/@{-} "a1";"b2",\ar@/_3pt/@{-} "b1";"a2"
    \vtwistneg~{"a2"}{"b2"}{"a3"}{"b3"}\vtwistneg~{"a3"}{"b3"}{"a4"}{"b4"}
    \vtwistneg~{"a4"}{"b4"}{"a5"}{"b5"}\vtwistneg~{"a5"}{"b5"}{"a6"}{"b6"}
    \xoverv~{"a7"}{"a6"}{"a8"}{"c7"}\xoverv~{"b8"}{"c7"}{"b7"}{"b6"}
    \hcap~{"b1"}{"b1"+(3,0.05)}{"b7"}{"b7"+(2,-.05)}
    \hcap~{"a1"}{"a1"-(3,-0.05)}{"a7"}{"a7"-(2,0.05)}|>
    \vcap~{"a8"-(0,1)}{"b8"-(0,1)}{"a8"}{"b8"}
\endxy
\quad
\begin{vmatrix}
    -1&0&-t(st-1)&s^2t^2-st+1\\
    0&-1&st&-s(st-1)\\
    -1&2st-2&-2st+3&0\\
    0&2st-1&-2st+2&-1
\end{vmatrix}
$$
\caption{\label{Fi:VKPairs} Virtual knot pairs with identical Alexander biquandles.}
\end{figure}

\section{The Biquandle of the Kishino Knot} \label{C:BQKishino}
We leave the reader with some open questions which remain.

\begin{figure}
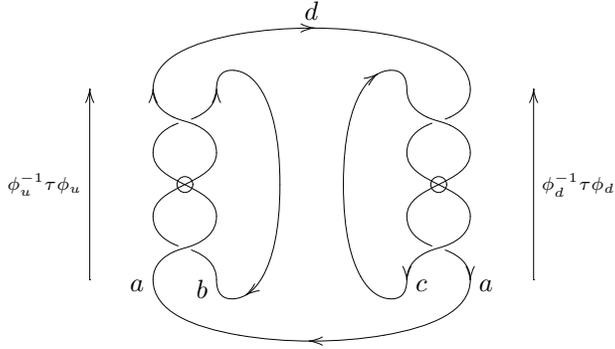

\[
{\xy 0;/r2pc/:
 (-2.5,0)="a1" *!UR{a\ }, (-2.5,1)="a2", (-2.5,2)="a3", (-2.5,3)="a4"
,(-1.5,0)="b1" *!UR{b\ }, (-1.5,1)="b2", (-1.5,2)="b3", (-1.5,3)="b4"
,(1.5,0)="c1" *!UL{\ c}, (1.5,1)="c2", (1.5,2)="c3", (1.5,3)="c4"
,(2.5,0)="d1" *!UL{\ a}, (2.5,1)="d2", (2.5,2)="d3", (2.5,3)="d4"
,{\save "b2"="temp","a2";"a3":"temp"::(0,0);(1,1)
     **\crv{(0.5,0)&(0.5,1)},(0,1);(1,0)
     **\crv{(0.5,1)&(0.5,0)},(0.5,0.5)
     *{\cir<3pt>{}}
  \restore}
,{\save "d2"="temp","c2";"c3":"temp"::(0,0);(1,1)
     **\crv{(0.5,0)&(0.5,1)},(0,1);(1,0)
     **\crv{(0.5,1)&(0.5,0)},(0.5,0.5)
     *{\cir<3pt>{}}
  \restore}
\vtwist~{"a2"}{"b2"}{"a1"}{"b1"}
\vtwist~{"c2"}{"d2"}{"c1"}{"d1"}=>
\vtwistneg~{"a4"}{"b4"}{"a3"}{"b3"}=<
\vtwistneg~{"c4"}{"d4"}{"c3"}{"d3"}
   \vcap~{"a4"+(0,2)}{"d4"+(0,2)}{"a4"}{"d4"}|< \POS?(.5)*+!D{d}
   \vcap~{"d1"-(0,2)}{"a1"-(0,2)}{"d1"}{"a1"}|<
   \vcap~{"b4"+(0,1)}{(-0.5,4.5)}{"b4"}{(-0.5,1.5)}
   \vcap~{(-0.5,-1.5)}{"b1"-(0,1)}{(-0.5,1.5)}{"b1"}|<
   \vcap~{(0.5,4.5)}{"c4"+(0,1)}{(0.5,1.5)}{"c4"}|<
   \vcap~{"c1"-(0,1)}{(0.5,-1.5)}{"c1"}{(0.5,1.5)}
,\ar^{\phi_u^{-1}\tau\phi_u} "a1"-(1,0);"a4"-(1,0)
,\ar_{\phi_d^{-1}\tau\phi_d} "d1"+(1,0);"d4"+(1,0)
\endxy}
\]
\caption{A non-trivial virtual knot (Kishino's example)}
\label{Fi:Kishino}
\end{figure}

There is an interesting knot discovered by Kishino
(see~\ref{Fi:Kishino}) which is nontrivial, but many attempts using
known virtual knot invariants fail to detect it.  In fact, it has
been shown to be nontrivial by a variety of methods ~\cite{FennBart,
DyeKauffman, KishinoSatoh, ManturovKauffman}. One of the more
intriguing methods is via the use of the quaternionic biquandle,
discovered by Bartholomew and Fenn ~\cite{FennBart}.

Here we give a version of that technique, by first obtaining a
presentation of the biquandle of the Kishino knot using our methods
and then showing that the resulting biquandle has a non-trivial
quaternionic representation.

Using biquandle morphisms, any biquandle associated with Kishino's
knot can be described with 3 generators, along with the following
relations which are computed from the diagram in~\ref{Fi:Kishino}:
\begin{align}
\phi_u^{-1}\tau\phi_u(a,b)&=(d,b) \\
\phi_d^{-1}\tau\phi_d(c,a)&=(c,d) \\
\end{align}
Calculating further,
\begin{align*}
\phi_u^{-1}\tau\phi_u\ (a,b)&=\phi_u^{-1}\tau\ (\ b\ur a\ ,\ a\lr b\ ) \\
&=\phi_u^{-1}\ (\ a\lr b\ ,\ b\ur a\ ) \\
&=(\ b\ur a\Bll{a\lr b}\ ,\ a\lr b\Bul{b\ur a}\ ) \\
\intertext{and}
\phi_d^{-1}\tau\phi_d\ (c,a)&=\phi_d^{-1}\tau\ (\ a\ul c\ ,\ c\ll a\ ) \\
&=\phi_d\ (\ c\ll a\ ,\ a\ul c\ ) \\
&=(a\ul c\Blr{c\ll a},\ c\ll a\Bur{a\ul c}\ ) \\
\end{align*}
so that the relations for the biquandle of the virtual knot diagram
in~\ref{Fi:Kishino} are:
\begin{align*}
b &= a\lr b\Bul{b\ur a} \\
c &= a\ul c\Blr{c\ll a} \\
d &= b\ur a\Bll{a\lr b} = c\ll a\Bur{a\ul c} \\
\end{align*}
This gives an associated biquandle via a description which has
3 generators and 3 relations:
\begin{equation} \label{E:KBQ}
\bracket{a,b,c\ |\ b = a\lr b\Bul{b\ur a},\ c = a\ul c\Blr{a\ll c},\ b\ur a\Bll{a\lr b} = c\ll a\Bur{a\ul c}}
\end{equation}
This description for a biquandle associated with Kishino's knot is
much simpler than the description obtained from the original
definition in Section~\ref{C:Biquandle}.
\begin{thm}
The description of the biquandle associated with Kishino's knot is non-trivial.
\end{thm}

The calculations in our proof are simplified by starting with the
description in \eqref{E:KBQ}. This demonstrates the efficiency the morphism
approach given in this paper.

The quaternionic biquandle is defined by the following rules:
\begin{align*}
a\ur b &= i a + (i+j) b \\
a\ul b &= i a + (1-j) b \\
a\lr b &= -i a + (i+j) b \\
a\ll b &=-i a + (1-j) b \\
\end{align*}

Applying this representation to the biquandle description in \eqref{E:KBQ},
the 3 equations become the following relations:
\begin{align}
(-3a+b) - 2 i b - 2 k b &= 0 \\
-(a+c)+ i (a+c) + k (a+c) &= 0 \\
(3b-3c)-4 i a &= 0 \\
\end{align}

Observe that by tensoring with the integers modulo three, denoted by $\Z_3$,
we obtain a module over the group ring of the quaternion group with mod 3 coefficients.
Note that by direct calculation, this is a non-trivial module,
for by direct reduction, we see that the equations above become

\begin{align*}
b - 2 i b - 2 k b &= 0 \\
-(a+c)+ i (a+c) + k (a+c) &= 0 \\
- i a &=0. \\
\end{align*}
Hence, $a=0$, giving us
\begin{align*}
(1-2i-2k)b &= 0 \\
(-1+i+k)c &= 0. \\
\end{align*}
This module is clearly non-trivial.  Hence, the Kishino knot is
non-trivial.

\bibliographystyle{abbrv}
\bibliography{BQPaper}

\end{document}